\documentclass[11pt]{article}

\usepackage{amssymb,latexsym,amsmath}
\usepackage{xcolor}
\usepackage{bm}

\usepackage{hyperref}

\begin{document}

\newtheorem{theorem}{Theorem}[section]

\newtheorem{proposition}[theorem]{Proposition}

\newtheorem{lemma}[theorem]{Lemma}

\newtheorem{corollary}[theorem]{Corollary}

\newtheorem{definition}[theorem]{Definition}

\newtheorem{remark}[theorem]{Remark}

\newtheorem{exempl}{Example}[section]

\newenvironment{exemplu}{\begin{exempl}  \em}{\hfill $\surd$

\end{exempl}}

\newcommand{\ea}{\mbox{{\bf a}}}
\newcommand{\eu}{\mbox{{\bf u}}}
\newcommand{\ep}{\mbox{{\bf p}}}
\newcommand{\ed}{\mbox{{\bf d}}}
\newcommand{\eD}{\mbox{{\bf D}}}
\newcommand{\eK}{\mathbb{K}}
\newcommand{\eL}{\mathbb{L}}
\newcommand{\eB}{\mathbb{B}}
\newcommand{\ueu}{\underline{\eu}}
\newcommand{\ueo}{\overline{u}}
\newcommand{\oeu}{\overline{\eu}}
\newcommand{\ew}{\mbox{{\bf w}}}
\newcommand{\ef}{\mbox{{\bf f}}}
\newcommand{\eF}{\mbox{{\bf F}}}
\newcommand{\eC}{\mbox{{\bf C}}}
\newcommand{\en}{\mbox{{\bf n}}}
\newcommand{\eT}{\mbox{{\bf T}}}
\newcommand{\eV}{\mbox{{\bf V}}}
\newcommand{\eU}{\mbox{{\bf U}}}
\newcommand{\ev}{\mbox{{\bf v}}}
\newcommand{\eve}{\mbox{{\bf e}}}
\newcommand{\uev}{\underline{\ev}}
\newcommand{\eY}{\mbox{{\bf Y}}}
\newcommand{\eP}{\mbox{{\bf P}}}
\newcommand{\eS}{\mbox{{\bf S}}}
\newcommand{\eJ}{\mbox{{\bf J}}}
\newcommand{\leb}{{\cal L}^{n}}
\newcommand{\eI}{{\cal I}}
\newcommand{\eE}{{\cal E}}
\newcommand{\hen}{{\cal H}^{n-1}}
\newcommand{\eBV}{\mbox{{\bf BV}}}
\newcommand{\eA}{\mbox{{\bf A}}}
\newcommand{\eSBV}{\mbox{{\bf SBV}}}
\newcommand{\eBD}{\mbox{{\bf BD}}}
\newcommand{\eSBD}{\mbox{{\bf SBD}}}
\newcommand{\ecs}{\mbox{{\bf X}}}
\newcommand{\eg}{\mbox{{\bf g}}}
\newcommand{\paromega}{\partial \Omega}
\newcommand{\gau}{\Gamma_{u}}
\newcommand{\gaf}{\Gamma_{f}}
\newcommand{\sig}{{\bf \sigma}}
\newcommand{\gac}{\Gamma_{\mbox{{\bf c}}}}
\newcommand{\deu}{\dot{\eu}}
\newcommand{\dueu}{\underline{\deu}}
\newcommand{\dev}{\dot{\ev}}
\newcommand{\duev}{\underline{\dev}}
\newcommand{\weak}{\rightharpoonup}
\newcommand{\weakdown}{\rightharpoondown}
\renewcommand{\contentsname}{ }

\title{Comments on Symplectic bipotentials}
\author{\href{https://imar.ro/~mbuliga/index.html}{Marius Buliga} \\ "Simion Stoilow" Institute of Mathematics of the Romanian Academy}
\date{Version: 18.04.2026, updated from version 16.02.2026}


\maketitle

\begin{abstract}
This is a reaction to the article Symplectic bipotentials, in published form \cite{hbgdspub} Harakeh M, Ban M, de Saxce G. Symplectic bipotentials. Mathematics and Mechanics of Solids. 2026;0(0) doi:10.1177/10812865251413554 , and in preprint form \cite{hbgds} arXiv:2410.23122v1. 

We give evidence that most of the content of the article \cite{hbgdspub} is already covered in previous works, partially cited  like \cite{bham} arXiv:0810.1419 [math.FA], or uncited, like  \cite{sben11} arXiv:1902.04598 [math-ph], \cite{bdiss} arXiv:2304.14158 [math-ph], which already introduced and studied symplectic bipotentials. 

These comments also apply to the conference paper version of  \cite{hbgds} arXiv:2410.23122v1, namely to the article \cite{hbgdspub2} Harakeh, M., Ban, M., de Saxce, G. (2026). Symplectic Bipotentials for the Dynamics of Dissipative Systems with Non Associated Constitutive Laws. In: Nielsen, F., Barbaresco, F. (eds) Geometric Science of Information. GSI 2025. Lecture Notes in Computer Science, vol 16034. Springer, Cham. doi:10.1007/978-3-032-03921-7\_31 .

\end{abstract}

\section{Introduction}

Most of the content of the article \cite{hbgdspub} is already covered in previous works:
\begin{enumerate}
\item[-] either cited, but presented in a misleading way, for example \cite{bham}
\item[-] or uncited, for example \cite{sben11} and \cite{bdiss}.
\end{enumerate}

Symplectic bipotentials, which  make the name of this article, are not new. They were already introduced in the uncited \cite{sben11}    then in the uncited \cite{bdiss} Definition 2.5 and studied in Theorem 2.6.

Section 5 BEN principle for symplectic bipotentials, in particular the SBEN principle, relation (24), is not new. It is already covered in the uncited \cite{bdiss}.

The symplectic Fenchel polar is not new. It is already introduced in \cite{MBGDS1}. 

Detailed comments are further presented in section \ref{detailedcomments}. 

In the following sections I cite as much as necessary from the mentioned articles, in order to help the reader to verify the claims.

\paragraph{Context.} Immediately after the appearance of the preprint version \cite{hbgds} arXiv:2410.23122v1 of the article under discussion, I posted a public commentary \cite{publiccom1} at my open notebook. The hope was that the authors either revise or retract the preprint version. I prepared a first version, in article form, of the same comments, and the present one is only a revised version of that. 

Recently, without knowing about the publication \cite{hbgdspub} of the article under discussion, I tried to see if my claims made in \cite{publiccom1} can be automatically verified. The experiment was successful in my opinion and the result was the second public commentary \cite{publiccom2}, which is largely automatically generated, with corrections clearly marked [between brackets]. The verification \cite{publiccom2} is not necessary for the content of these comments, it only shows  that clear formulated claims can, with some effort, be checked by present available versions of LLMs, therefore these claims can also be checked by willing independent readers.

As a matter of principle, the idea is to not believe any claim by default, instead to provide as many means as possible to the reader in order to make an independent verification. 

These comments were sent to the Editor of the journal Mathematics and Mechanics of Solids (MMS), where the published version \cite{hbgdspub} appeared. 

\paragraph{Update 1.} After the appearance of the first version of these Comments was discovered the conference paper version of  \cite{hbgds} arXiv:2410.23122v1, namely  the article \cite{hbgdspub2}. This paper is more informal than \cite{hbgdspub} but it makes the same originality claims which are discussed in these Comments.

\paragraph{Update 2.}  Following the first version of these Comments, the Chief Editor of the journal Mathematics and Mechanics of Solids sent a letter of apology and informed me this letter will be included in the next available issue of MMS. As there is lag of about a year between online first and published issues of MMS, it is not clear when this leter will be published, while in the same time the article  \cite{hbgdspub} is presented in MMS (online first) as original research.  That is why the letter of the editor is reproduced here (\href{https://chorasimilarity.wordpress.com/wp-content/uploads/2026/03/from-the-editor.pdf}{link to pdf}). 

\vspace{.5cm}

"From the Editor: \\

After the online publication of [1] I received a complaint from M. Buliga [2] in which he claimed that certain aspects of its contents had been discussed in his prior work, which had not been cited. Upon reviewing the matter, I have concluded that Buliga’s concern is warranted. If I had been aware of this issue, I would have asked the authors of [1] to include citation and discussion of Buliga’s relevant contributions in their paper. I apologize for this oversight.

\vspace{.2cm}

David Steigmann, Editor-in-Chief \\

\vspace{.2cm}

[1] M. Harakeh, M. Ban and G. de Saxce, 2026, Symplectic bipotentials. Math. Mech. Solids (\href{https://doi.org/10.1177/10812865251413554}{doi:10.1177/10812865251413554})

[2] M. Buliga, Comments on symplectic bipotentials. \href{https://arxiv.org/abs/2602.14614}{https://arxiv.org/abs/2602.14614} "

\vspace{.5cm}

In a further email exchange. to my inquiry about when this letter will be published and why the editor keeps the article as original reseach, the editor informed me that: "As I stated, the letter will be published as soon as possible. The matter will then be closed as far as MMS is concerned. Any reasonable person would conclude by reading it that the paper in question is not original work."

\paragraph{Acknowledgements.}The present version of the comments takes into consideration the published version of the article under discussion. I want to thank the Chief Editor of MMS for sending me the pdf of the published article. There are almost no differences between the arXiv version and the published version and all comments apply to both \cite{hbgds} and \cite{hbgdspub}, as well as to the conference version \cite{hbgdspub2}.

\section{Preliminaries}

Hamiltonian inclusions with convex dissipation provide an extension of hamiltonian mechanics to irreversible processes \cite{bham},  \cite{sben1}, \cite{sben11}, \cite{bdiss}. In \cite{MBGDS1} the alternative name "SBEN" (ie symplectic Brezis-Ekeland-Nayroles) was used, because, in some applications, the quasistatic formal approximation of hamiltonian inclusions allowed us to recover the old variational principles of Brezis-Ekeland \cite{Brezis Ekeland 1976} and Nayroles \cite{Nayroles 1976}. However, in general the quasistatic approximation of hamiltonian inclusions  contradicts the mathematical structure of the formulation, which is naturally dynamical.

In hamiltonian mechanics a physical system is described by a state vector $q\in Q$ and a momentum vector $p \in P$. The evolution in time of the system is governed by a hamiltonian function $H=H(q,p,t)$, via the equations: 
\begin{equation} \left\{   
\begin{array}{rcl}
\dot{q} & = & \frac{\partial H}{\partial p}  (q,p,t) \\
- \dot{p} & = &  \frac{\partial H}{\partial q}  (q,p,t)
\end{array}
\right.
\label{hamilton}
\end{equation} 
where $\displaystyle \dot{q}$, $\displaystyle \dot{p}$ denote derivatives with respect to time. The evolution  is reversible. 

The symplectic gradient of $H$, denoted by $\displaystyle XH(q,p,t) \in Q \times P$,  is 
\begin{equation}
XH(q,p,t) = \left(\frac{\partial H}{\partial p}  (q,p,t), -\frac{\partial H}{\partial q}  (q,p,t) \right)
\label{sympgrad}
\end{equation}

The evolution equation (\ref{hamilton}) for hamiltonian dynamics is therefore: 
$$ \dot{c}(t) \, = \, XH(c(t),t)$$ 

A dissipative modification of hamiltonian mechanics, called hamiltonian inclusions with convex dissipation,   \cite{sben11}, \cite{bdiss}, which extends the work from \cite{bham}, \cite{MBGDS1}, \cite{sben1},  is obtained by the introduction of the gap vector  $\displaystyle \eta = (\eta_{q}, \eta_{p})$. We modify the hamiltonian dynamics equations (\ref{hamilton}) in the following way: 

\begin{equation} \left\{   
\begin{array}{rcl}
\dot{q} & = & \frac{\partial H}{\partial p}  (q,p,t) \, + \, \eta_{q} \\
\dot{p} & = &  - \frac{\partial H}{\partial q}  (q,p,t) \, + \, \eta_{p} 
\end{array}
\right.
\label{hamiltongap}
\end{equation} 
The equations for the gap vector are introduced via likelihoods \cite{bdiss} section 2, here Definition \ref{deflikelihood}. In the symplectic space $\displaystyle Q \times P$ we introduce the maximal likelihood between two vectors as: 
\begin{equation}
\pi_{max}(z',z") \, = \, e^{\min \left\{0, \omega(z',z") \right\}}
\label{maxlike}
\end{equation} 
This is a number in $(0,1]$, for example when $\displaystyle z'$ and $\displaystyle z"$ are colinear their maximal likelihood is 1, or if $\displaystyle  \omega(z',z") > 0$ then again the maximal likelihood is 1. 

\begin{definition}
A likelihood function is a function $\displaystyle \pi: (Q \times P)^{3} \rightarrow [0,1]$ with the properties: for any $\displaystyle z,z',z" \in Q \times P$
\begin{enumerate}
\item[(a)] if either of the  maxima exist: $\displaystyle \max_{v \in Q \times P} \pi(z, z', v)  \, , \, \max_{w \in Q \times P} \pi(z, w, z")$,  then they are equal to $0$ or $1$
\item[(b)] the functions 
$\displaystyle v \in Q \times P \mapsto - \ln \pi(z, z', v) \, , \, w \in Q \times P \mapsto - \ln \pi(z, w, z")$
are convex and lower semicontinuous (lsc).
\end{enumerate} 

A likelihood $\displaystyle \pi: (Q \times P)^{3} \rightarrow [0,1]$ is tempered if moreover for any $\displaystyle z,z',z" \in Q \times P$
\begin{enumerate}
\item[(c)] $ \displaystyle \pi(z,z',z") \leq \pi_{max}(z',z")$
\end{enumerate}
\label{deflikelihood}
\end{definition} 

We define the information content function associated to the likelihood $\pi$ as  
\begin{equation}
I: (X \times Y)^{3} \rightarrow [0,+\infty] \, \, , \, \, I(z,z',z") \, = \, - \ln \pi(z, z', z")
\label{icfun}
\end{equation}
with the convention that $\displaystyle - \ln 0 \, = \, + \infty$.

With likelihoods, the dissipative modification is described further. 

\begin{definition}
Given a hamiltonian $H$ and a tempered likelihood function $\pi$ 
$$\displaystyle H: Q \times P \times \mathbb{R} \rightarrow \mathbb{R} \quad , \quad \pi: (Q \times P)^{3} \rightarrow [0,1]$$ 
the dynamics of a  physical system is defined by the modification of the hamiltonian dynamics equations (\ref{hamiltongap}), 
together with the new equation:  
\begin{equation}
\pi(z, \dot{z}, \eta) \, = \, 1
\label{likely}
\end{equation} 
with the notations $\displaystyle z=(q,p), \eta = (\eta_{q}, \eta_{p}) \in Q \times P$, $\displaystyle \dot{z} = (\dot{p}, \dot{q}) \in Q \times P$. 
\label{mainproblem}
\end{definition}

In \cite{sben11}  Proposition 1.3 (here \ref{biposync}) is introduced  the function 
\begin{equation}
b: (X \times Y)^{3} \rightarrow \mathbb{R} \cup \left\{ +\infty \right\} \, \, , \, \, b(z,z',z") \, = \ I(z,z', z") \, + \, \omega(z', z")
\label{bipoics}
\end{equation} 

The gap vector $\displaystyle \eta$ satisfies therefore  a new equation: 
\begin{equation}
b(z,\dot{z},\eta) = \omega(\dot{z},\eta)
\label{gapequation}
\end{equation}
which is formulated in terms of a function $\displaystyle 
b: (Q \times P)^{3} \rightarrow \mathbb{R} \cup \left\{ +\infty \right\}$ 
which {\bf has to be} a bipotential (definition \ref{biposync}) in its 2nd and 3rd variables, i.e. it satisfies: for any $z \in Q \times P$,  
\begin{enumerate}
\item[(a)] for any $z', z" \in Q \times P$ the functions $b(z, z', \cdot)$ and $b(z, \cdot, z")$  are  convex (and lsc), 
\item[(b)] for any $z', z" \in Q \times P$ we have the equivalences 
\begin{equation}
z' \in \, \partial b(z, z', \cdot) (z") \, \Longleftrightarrow \, z" \in \, \partial b(z, \cdot, z") (z') \, \Longleftrightarrow \, b(z,z',z") = \langle\langle z',z"\rangle\rangle
\label{bipoequiv2}
\end{equation}
\end{enumerate}
Moreover the function $b$ is non-negative: $\displaystyle b(z,z',z") \geq 0$ for any $\displaystyle z, z', z" \in Z$ because it comes from a tempered likelihood. This is important for  
the dissipation inequality, see \cite{bdiss} Definition 3.3. (here \ref{here33}) and Theorem 3.4 (here \ref{theorembalance}).

\section{Detailed comments}
\label{detailedcomments}

Most of the content of the article \cite{hbgdspub} is already covered in previous works:
\begin{enumerate}
\item[-] either cited, but presented in a misleading way, for example \cite{bham}
\item[-] or uncited, for example \cite{sben11} and \cite{bdiss}.
\end{enumerate}

Symplectic bipotentials, which  make the name of this article, are not new. They were already introduced in the uncited \cite{sben11}    then in the uncited \cite{bdiss} Definition 2.5 and studied in Theorem 2.6.

Section 5 BEN principle for symplectic bipotentials, in particular the SBEN principle, relation (24), is not new. It is already covered in the uncited \cite{bdiss}.

The symplectic Fenchel polar is not new. It is already introduced in \cite{MBGDS1}.

The article \cite{bham}   Definition  2.2. (here \ref{defssub}) introduces for the first time the symplectic subdifferential.

\begin{definition}
Let $F: X \times Y \rightarrow \mathbb{R}$ be a convex lsc function. The
symplectic subdifferential of $F$ is the multivalued function which sends 
$z = (x,y) \in X \times Y$ to the set 
$$\displaystyle X \, F (z) \, = \, \left\{ 
z' \in X \times Y \mbox{ : } \forall \, z" \in X \times Y \quad  
F(z+z") \,  \geq \,  F(z) \, + \, 
\omega (z' , z")\right\} $$ 
\label{defssub}
\end{definition}

In the article \cite{hbgdspub} the reference is generically  cited at the beginning of Section 2.1, along with a much more recent article \cite{MBGDS1} with de Saxce. In relation (2), where the symplectic subdifferential is defined, there is no mention of the source.

Moreover, in the following relation (3), an element of the symplectic subdifferential is described to satisfy the “so-called Hamiltonian inclusion”. The relation (3) is not a Hamiltonian inclusion.  Hamiltonian inclusions are described in the introduction, in relation (1),  which is misleading, compare with  \cite{bham} Definition 2.3 (here \ref{mainevol}) and relation (14), here (\ref{mainsymp}).

\begin{definition}
Let $H : [0,T] \times X \times Y \rightarrow \mathbb{R}$ such that for all 
$t \in [0,T]$ we have $H(t, \cdot) \in Der(X,Y)$, and $\displaystyle 
D: (X \times Y)^{2} \rightarrow \mathbb{R} \cup \left\{ + \infty \right\}$ be a
"dissipation function" with the properties: 
\begin{enumerate}
\item[(a)] for any $z', z" \in X \times Y$  we have $\mathcal{R}(z', z") \geq 0$ and 
$\mathcal{R}(z' ,0) = 0$, 
\item[(b)] for any $z \in X \times Y$ the function $\mathcal{R}(z,  \cdot)$ is convex, 
lsc. 
\end{enumerate}
Then a curve $z: [0,T] \rightarrow X \times Y$ is a solution of the 
evolution problem with Hamiltonian
$H$ and dissipation $D$ if it is derivable for all $t \in [0, T]$ (with
differential denoted by $\displaystyle \dot{z}$) and it
satisfies the subdifferential inclusion: 
\begin{equation}
\dot{z}(t) \, - \, X_{H(t,\cdot)} (z(t)) \, \in \, X\, \left( \mathcal{R}(z(t), \cdot )
\right) (\dot{z}(t)) \quad \quad .  
\label{mainsymp}
\end{equation}
\label{mainevol}
\end{definition}

A Hamiltonian inclusion is expressed as the irreversible part of the time rate of the evolution belongs to the symplectic subgradient of the dissipation, while the reversible part of the time rate equals the symplectic gradient of the Hamiltonian.

After the  relation (1)  the authors write about \cite{bham} “Nevertheless, it is important to remark that [the dissipation] is not necessarily 1-homogeneous or even homogeneous”. But later they write “To release the restrictive hypothesis of 1-homogeneity (in particular to address viscoplasticity), we introduce in this
work the symplectic Fenchel polar.”

The 1-homogeneity of the dissipation is a hypothesis in Mielke theory of quasistatic evolutionary processes. Again in \cite{bham} Mielke quasistatic theory is compared with the dynamical Hamiltonian inclusions. In the uncited  \cite{bdiss} “Relation with rate-independent processes” this is again explained with details (see further more about this uncited paper).

There is no limitation of 1-homogeneity in the dynamical theory! This is not why symplectic Fenchel polar is introduced.

By the way, this symplectic Fenchel polar is already introduced in \cite{MBGDS1}, again in Definition 2.2 (here \ref{dspolar}), so this is not introduced in the present work.

\begin{definition}
The symplectic polar, or the symplectic Fenchel transform of $F: X \times Y \rightarrow \mathbb{R}\cup \left\{+\infty\right\}$, a convex lsc function, is the function: 
$$F^{*\omega}(z') \, = \, \sup \left\{ \omega(z',z) - F(z) \mbox{ : } z \in X \times Y \right\}$$
\label{dspolar}
\end{definition}

A short comment on \cite{MBGDS1} is necessary. This article re-tells the story of Hamiltonian inclusions from \cite{bham} and then makes the original observation that when we pass from the dynamical to the quasistatic regime, we recover the Brezis-Ekeland and Nayroles principles. It is misleading to attribute the Hamiltonian inclusions to this article, when the only novelty is the connection with the BEN principle.

Symplectic bipotentials  are defined in this paper in Section 4 Symplectic bipotentials p. 7 (a) (b) (c).

But they were already introduced in the uncited \cite{sben11}   where they appear first in Proposition 1.3 (here \ref{biposync}), after the definition of  the information content (4), here (\ref{icfun}) and of the required properties from  Definition 1.2 (here \ref{likelyprop}).

\begin{definition}
The likelihood $\pi$ and the associated information content $I$ satisfy: 
\begin{enumerate}
\item[(a)] The information content function (\ref{icfun}) is convex in each of the 2nd and 3rd variables and it has the needed degree of smoothness required (for example it is lower semi-continuous with respect to the relevant topologies on $X$ and $Y$)  
\item[(b)] for any $\displaystyle z, z' \in X \times Y$,  the following maxima exist 
$$\displaystyle \max_{z" \in X \times Y} \pi(z, z', z")  \, , \, \max_{z" \in X \times Y} \pi(z, z", z')$$ and they are either $0$ or $1$.
\end{enumerate}
\label{likelyprop}
\end{definition}

\begin{proposition}
The information content function (\ref{icfun}) satisfies the conditions  from Definition \ref{likelyprop}  if and only if the function  
\begin{equation}
b: (X \times Y)^{3} \rightarrow \mathbb{R} \cup \left\{ +\infty \right\} \, \, , \, \, b(z,z',z") \, = \ I(z,z', z") \, + \, \omega(z', z") 
\label{bipoic}
\end{equation} 
is a bipotential, i.e. it satisfies: for any $z \in X \times Y$,  
\begin{enumerate}
\item[(a)] for any $z', z" \in X \times Y$ the functions $b(z, z', \cdot)$ and $b(z, \cdot, z")$  are  convex (and lsc), 
\item[(b)] for any $z', z" \in X \times Y$ we have the equivalences 
\begin{equation}
z' \in \, \partial b(z, z', \cdot) (z") \, \Longleftrightarrow \, z" \in \, \partial b(z, \cdot, z") (z') \, \Longleftrightarrow \, I(z, z', z") = 0
\label{bipoequiv}
\end{equation}
where "$\partial$" denotes a subgradient. 
\end{enumerate}
\label{biposync}
\end{proposition}

Then in the uncited \cite{bdiss} Definition 2.5 (here \ref{bipodef})and studied in Theorem 2.6 (here \ref{likebipo}).  

For a space $\displaystyle Q \times P$ which is in duality $\displaystyle d: \left(Q \times P\right)^{2} \rightarrow \mathbb{R}$ with itself, bipotentials have the following definition. 

\begin{definition}
A function $\displaystyle b: \left(Q \times P\right)^{2} \rightarrow \mathbb{R} \cup \left\{+\infty\right\}$ is a bipotential if: 
\begin{enumerate}
\item[(a)] $\displaystyle b(z',z") \geq d(z',z")$ for any $\displaystyle z',z" \in Q \times P$
\item[(b)] $\displaystyle b(z',z") =  d(z',z")$ if and only if $\displaystyle z" \in \partial_{d}^{R} b(\cdot,z")(z')$ if and only if $\displaystyle z' \in \partial_{d}^{L} b(z',\cdot)(z")$
\end{enumerate}
\label{bipodef}
\end{definition}

Here appear left $\displaystyle \partial_{d}^{L}$ and right $\displaystyle \partial_{d}^{R}$ subgradients, which are defined with respect to a general duality $\displaystyle d$ in the uncited \cite{bdiss} Definition 2.2, and then in Theorem 2.3 (Fenchel inequalities) there is a general formulation of the Fenchel inequalities in this context. This is to be compared with the article under discussion, section 2, where a more informal presentation of the same ideas is done for the particular case when the duality is a symplectic form, thus antisymmetric. 

Likelihoods are related to bipotentials. The relation has been noted before, where information contents of likelihoods appear as syncs, definition 2.3 \cite{gds4}. The same was first observed in relations (51), (52) \cite{laborderenard}. The proof of the following theorem is the same as the one of proposition 2.4 \cite{gds4}, only adapted to the notations of the present article. 

\begin{theorem}
Let $\displaystyle \pi: (Q \times P)^{3} \rightarrow [0,1]$ be a function and $\displaystyle d: \left(Q \times P\right)^{2} \rightarrow \mathbb{R}$ be a duality of $\displaystyle Q \times P$ with itself. Denote by $\displaystyle I = - \ln \pi$ the information content of the function $\displaystyle \pi$, and by 
$$b_{d}(z,z',z") = I(z,z',z") + d(z',z")$$
The function $\displaystyle \pi$ is a likelihood if and only if the function $\displaystyle b_{d}(z,\cdot,\cdot)$ is a bipotential with respect to the duality $\displaystyle d$.  
\label{likebipo}
\end{theorem}

As we can see, while likelihoods are independent of dualities, bipotentials are relative to a duality. We can easily transform a bipotential $\displaystyle b$ with respect to the duality $\displaystyle d$ into another bipotential $\displaystyle b'$ with respect to another duality $d'$, by the formula: 
$$b'(z',z") - d'(z',z") \, = \, b(z',z") - d(z',z")$$

For the particular duality $\displaystyle \omega$, the symplectic form, to any likelihood we associate its symplectic bipotential. 

\begin{definition}
Let $\displaystyle \pi: (Q \times P)^{3} \rightarrow [0,1]$ be a likelihood. 
With $\displaystyle I(z,z',z") = - \ln \pi(z,z',z")$ the information content of $\displaystyle \pi$, the symplectic bipotential associated to this likelihood is
$$b_{\omega}^{\pi}(z,z',z") \, = \, I(z,z',z") + \omega(z',z")$$

The minimal symplectic bipotential is the symplectic bipotential of the maximal likelihood $\displaystyle \pi_{max}$, i.e.
$$b_{\omega}^{min}(z',z") \, = \,  \max \left\{0, \omega(z',z") \right\}$$
\label{bsimp}
\end{definition}

A likelihood $\displaystyle \pi$ is tempered if and only if for any $\displaystyle z, z',z" \in Q \times P$ 
$$b_{\omega}^{\pi}(z,z',z") \geq 0$$
or equivalently 
$$b_{\omega}^{\pi}(z,z',z") \geq b_{\omega}^{min}(z',z")$$

Section 5 BEN principle for symplectic bipotentials, in particular the SBEN principle, relation (24), is already covered in the uncited \cite{bdiss} Definition 3.3. (here \ref{here33}) and Theorem 3.4 (here \ref{theorembalance}).

Recall that for the same information content, for different dualities we obtain different bipotentials, therefore 
$$b_{d}(z,z',z") - d(z',z") \, = \, b_{\omega}^{\pi}(z,z',z") - \omega(z',z") \, = \, I(z,z',z")$$

\begin{definition}
The dissipation along a curve $\displaystyle t \in [0,T] \mapsto c(t) \in Q \times P$ is the functional: 
\begin{equation}
Diss^{\pi}(c,0,T) \, = \, \int_{0}^{T} b_{\omega}^{\pi}(\dot{c}(t), \dot{c}(t) - XH(c(t),t))\mbox{ d}t
\label{defdissip}
\end{equation}  
Remark that for any curve 
$$Diss^{\pi}(c,0,T) \, \geq \, 0$$
because the likelihood $\displaystyle \pi$ is tempered. 
\label{here33}
\end{definition}

In the following theorem we give the energy balance and  dissipation inequalities. 

\begin{theorem}
Let $\displaystyle t \in [0,T] \mapsto (c(t),\eta(t))$ be a solution of \ref{mainproblem} with the initial condition $\displaystyle c(0) = z_{0}$. Then for any $t \in [0,T]$: 
\begin{enumerate}
\item[(a)](energy balance) 
\begin{equation}
H(c(t),t) \, = \, H(z_{0},0) +  \int_{0}^{t} \frac{\partial H}{\partial t}(c(\tau),\tau)\mbox{ d}\tau \, - \, Diss^{\pi}(c,0,t)
\end{equation}
\item[(b)] (dissipation inequalities)
$$\frac{d}{dt} H(c(t),t)) -  \frac{\partial H}{\partial t}(c(t),t) \, \leq \, 0$$
For any curve $\displaystyle t \in [0,T] \mapsto c'(t)$ which satisfies $\displaystyle c'(0) = c(0)$ we have 
\begin{equation}
Diss^{\pi}(c',0,t) + H(c'(t),t) \, \geq Diss^{\pi}(c,0,t) + H(c(t),t)
\label{mindiss}
\end{equation}
\end{enumerate}
\label{theorembalance}
\end{theorem}

The Section  5.1. Application to the non associated plasticity, mirrors the application in  \cite{MBGDS1} to standard plasticity. And also the application to elastoplasticity from the uncited  \cite{sben11}  Section 2.  The only “new” thing is to use a symplectic bipotential which does not have to be associated, which is trivial in the light of the uncited papers.

\end{document}